\documentclass[12pt]{article}
\usepackage{amsmath}
\usepackage{amsthm}
\usepackage{amssymb}
\usepackage{epsfig}

\makeatletter
\newfont{\footsc}{cmcsc10 at 8truept}
\newfont{\footbf}{cmbx10 at 8truept}
\newfont{\footrm}{cmr10 at 10truept}
\makeatother 

\def\S #1{\mathcal{S}_{#1}}

\def\P #1{\mathcal{P}_{#1}}

\def\Qj {Q_j}
\def\T {\mathcal{T}}

\def\st {ST}
\def\stc {\mathcal{ST}}
\def\J {\mathcal{J}}

\def\itau {\tilde{\tau}}
\def\jtau {\overline{\tau}}

\newcommand{\siav}[2]{I^{s}_{#1}(#2)}
\newcommand{\opn}[1]{F({#1})}
\newcommand{\opi}[1]{G({#1})}

\theoremstyle{plain}
\newtheorem{theorem}{Theorem}[section]
\newtheorem{lemma}[theorem]{Lemma}
\newtheorem{prop}[theorem]{Proposition}
\newtheorem{cor}[theorem]{Corollary}

\theoremstyle{definition}
\newtheorem{defn}[theorem]{Definition}
\newtheorem*{q}{Question}

\newtheorem{ex}[theorem]{Example}

\theoremstyle{remark}
\newtheorem{remark}[theorem]{Remark}

\theoremstyle{plain}

\title{Subsequence containment by involutions}

\author{Aaron D.\ Jaggard\thanks{Part of this work is drawn from the
author's Ph.D. dissertation which was written under the direction
of Herbert S.\ Wilf at the University of Pennsylvania and
partially supported by ONR grant N00014--01--1--0431 and the DoD
University Research Initiative (URI) program administered by the
ONR under grant N00014--01--1--0795.  The part of this work
carried out at Tulane University was partially supported by NSF
grant DMS--0239996. The presentation of part of this work at the
Permutation Patterns 2003 conference was partially supported by
Penn GAPSA and
the New Zealand Institute for Mathematics and its Applications.}\\
\small Department of Mathematics\\[-0.8ex]
\small Tulane University\\[-0.8ex]
\small New Orleans, LA 70118 USA\\[-0.8ex]
\small \texttt{adj@math.tulane.edu}}

\date{\small Version of 26 August 2004\\
\small MR Subject Classifications: 05A05, 05A15, 05E10}

\begin{document}

\maketitle

\begin{abstract}
Inspired by work of McKay, Morse, and Wilf, we give an exact count
of the involutions in $\S{n}$ which contain a given permutation
$\tau\in\S{k}$ as a subsequence; this number depends on the
patterns of the first $j$ values of $\tau$ for $1\leq j\leq k$. We
then use this to define a partition of $\S{k}$, analogous to
Wilf-classes in the study of pattern avoidance, and examine
properties of this equivalence.  In the process, we show that a
permutation $\tau_1\ldots\tau_k$ is layered iff, for $1\leq j\leq
k$, the pattern of $\tau_1\ldots\tau_j$ is an involution.  We also
obtain a result of Sagan and Stanley counting the standard Young
tableaux of size $n$ which contain a fixed tableau of size $k$ as
a subtableau.
\end{abstract}

\section{Introduction}

Given a permutation $\pi=\pi_n\ldots\pi_n$ in the symmetric group
$\S{n}$ and a word $\sigma=\sigma_1\ldots\sigma_k$ of $k$ distinct
letters, say that $\pi$ \textit{contains $\sigma$ as a
subsequence} if $\pi_{i_1}=\sigma_1$, \ldots, $\pi_{i_k}=\sigma_k$
for some $i_1 < \cdots < i_k$.  Recent work of McKay, Morse, and
Wilf~\cite{McKay2000} implies that the probability that an
involution in $\S{n}$ contains any fixed $\tau\in\S{k}$ as a
subsequence is $1/k!+o(1)$ as $n\rightarrow\infty$.  We sharpen
this aspect of their work with the following theorem.
\newtheorem*{invthm}{Theorem~\ref{thm:inv}}
\begin{invthm}
For a fixed permutation $\tau=\tau_{1}\tau_{2}\ldots\tau_{k}
\in\S{k}$ and $n\geq k$, the number $\siav{n}{\tau}$ of
involutions in $\S{n}$ which contain $\tau$ as a subsequence
equals
\begin{equation*}
\sideset{}{'}\sum_{j}\binom{n-k}{k-j}t_{n-2k+j}
\end{equation*}
where the sum is taken over $j=0$ and those $j\in[k]$ such that
the pattern of $\tau_1\ldots\tau_j$ is an involution in
$\mathcal{S}_j$.
\end{invthm}
\noindent We use $[k]$ for $\{1,\ldots, k\}$; recall that the
\textit{pattern} of a word $\sigma=\sigma_1\ldots\sigma_k$ of $j$
distinct letters is the order-preserving relabeling of the letters
of $\sigma$ with $[j]$.  We will refer to the pattern of
$\sigma_1\ldots\sigma_j$ as the \textit{length $j$ initial pattern
of $\sigma$\/}.

We then define $\sigma$ and $\tau$ to be equivalent iff, for every
$n$, $\siav{n}{\sigma}=\siav{n}{\tau}$; this leads to an
apparently new classification of permutations.  Defining
$\J(\tau)$ to be the set of indices $j$ over which the sum in
Theorem~\ref{thm:inv} is taken, $\sigma$ and $\tau$ are equivalent
iff $\J(\sigma)=\J(\tau)$ and $\sigma,\tau\in\S{k}$ for some $k$.
We then examine some enumerative results relating to the sets
$\J(\tau)$ in general; this leads us to an apparently new
characterization of layered permutations.
\newtheorem*{layprop}{Proposition~\ref{prop:layered}}
\begin{layprop}
A permutation $\tau\in\S{k}$ is a layered permutation if and only
if $\J(\tau)=\{0,1,\ldots,k\}$, \emph{i.e.}, if and only if the
pattern of $\tau_1\ldots\tau_j$ is an involution for every $j\in
[k]$.
\end{layprop}
We also prove the following theorem about $\J(\tau)$ for those
$\tau$ corresponding to a given Young tableau.
\newtheorem*{ajlem}{Lemma~\ref{lem:aj}}
\begin{ajlem}
Fix $k$ and $\lambda\vdash k$, and let $T$ be a standard Young
tableau of shape $\lambda$.  Let $\T$ be the set of permutations
which correspond to $(T,Q)$ for some $Q$, let $a_{0}=f^{\lambda}$,
and for $1\leq j\leq k$ define $a_{j}$ to be the number of
$\tau\in\mathcal{T}$ for which $j\in\J(\tau)$.  Then (with
$f^{\emptyset}=1$)
\begin{equation*}
a_{j}=\sum_{\mu\,\vdash\, j}f^{\lambda/\mu} \quad (j=0,1,\ldots,
k),
\end{equation*}
where $f^{\lambda/\mu}$ is the number of standard Young tableaux
of skew shape $\lambda/\mu$.
\end{ajlem}
\noindent This leads to a result of Sagan and
Stanley~\cite{Sagan1990} counting the Young tableaux of size
$n\geq k$ in which the entries $1,\ldots, k$ form a specified
subtableau.

Section~\ref{sec:pqrand} contains Theorem~\ref{thm:inv} and
related material. Section~\ref{sec:app} covers enumerative
questions related to the sets $\J(\sigma)$ and connects these
ideas to previous work on Young tableaux.

\section{$P$-quasirandomness}
\label{sec:pqrand}

\subsection{Definition}
Quasirandom permutations were introduced by McKay, Morse, and
Wilf~\cite{McKay2000} and are defined as follows.
\begin{defn}
Let $\mathcal{P}_{n}\subseteq\S{n}$ be a non-empty set of
permutations for infinitely many values of $n$, and let
$\mathcal{P}=\cup_{n}\mathcal{P}_{n}$.  For a word $\sigma$ of $k$
distinct letters from $[n]$, let $h(n,\sigma)$ be the number of
permutations in $\mathcal{P}_{n}$ which contain $\sigma$ as a
subsequence.  If $\mathcal{P}_{n}\neq\emptyset$, define
$g(n,\sigma) = h(n,\sigma)/|\mathcal{P}_{n}|$, the probability
that $\pi\in\mathcal{P}_{n}$ contains $\sigma$ as a subsequence.
$\mathcal{P}$ is \textit{quasirandom} (or a \textit{quasirandom
family of permutations}) if, for every $k\geq 1$, we have
\begin{equation}
\lim_{n\rightarrow\infty}
\max_{\sigma}\left|g(n,\sigma)-\frac{1}{k!}\right|=0,
\label{eq:qrand}
\end{equation}
where the maximum is over all sequences $\sigma$ of $k$ distinct
elements of $[n]$ and the limit is over those $n$ such that
$\mathcal{P}_{n}\neq\emptyset$.
\end{defn}
\noindent McKay, Morse, and Wilf used the quasirandomness of the
set of involutions to prove theorems about the entries of Young
tableaux.  Their results only need Equation~\ref{eq:qrand} to hold
for $\sigma\in\S{k}$ and not necessarily for arbitrary words
$\sigma$ of $k$ distinct integers; this leads us to define the
strictly weaker notion of \textit{$p$-quasirandom}
(\textit{permutation-quasirandom}) permutations, whose definition
repeats that of quasirandom permutations except that the word
$\sigma$ is replaced by $\tau\in\S{k}$.
\begin{defn}
Let $\mathcal{P}_{n}\subseteq\S{n}$ be non-empty for infinitely
many values of $n$, and let $\mathcal{P}=\cup_{n}\mathcal{P}_{n}$.
For $\tau\in\S{k}$, let $h(n,\tau)$ be the number of permutations
in $\mathcal{P}_{n}$ which contain $\tau$ as a subsequence. Let
$f(n,\tau) = h(n,\tau)/|\mathcal{P}_{n}|$, the probability that a
permutation in $\P{n}$ contains $\tau$, if
$\mathcal{P}_{n}\neq\emptyset$. $\mathcal{P}$ is a
\emph{p-quasirandom} (\emph{permutation-quasirandom}) \emph{family
of permutations} if, for all $k\geq 1$
\begin{equation}
\lim_{n\rightarrow\infty}\max_{\tau}\left|f(n,\tau)-\frac{1}{k!}\right|=0,
\label{eq:pqrand}
\end{equation}
where the maximum is now taken over all $\tau\in\S{k}$ and the
limit is again restricted to be over values of $n$ for which
$\mathcal{P}_{n}\neq\emptyset$.
\end{defn}

It is clear that if $\mathcal{P}$ is quasirandom, then
$\mathcal{P}$ is $p$-quasirandom.  The following example shows
that the converse is not true, so $p$-quasirandomness is indeed
strictly weaker than quasirandomness.
\begin{ex}
Define
\begin{equation}
\mathcal{P}_{n}=
\begin{cases}
\{\pi\in\S{n}|\pi\mathrm{\ fixes\ }[n]\setminus[n/2]\}, &n\mathrm{\ even}\\
\emptyset, &n\mathrm{\ odd}
\end{cases}.\notag
\end{equation}
For fixed $k$, for every even $n\geq 2k$ we see that every
$\pi\in\mathcal{P}_{n}$ contains the $k$-element list
$\sigma_{n}=(n-k+1)(n-k+2)\cdots(n-1)n$ as a subsequence.  Thus
$g(n,\sigma_{n})=1$, and we have that the limit in
Equation~\ref{eq:qrand} is equal to $1-1/k!$ and
$\{\mathcal{P}_{n}\}$ is not a quasirandom family.

However, for fixed $k$ and every even $n\geq 2k$, the probability
that $\pi\in\mathcal{P}_{n}$ contains $\tau\in\S{k}$ as a
subsequence is just the probability that $\pi'$ chosen uniformly
at random from $\S{n/2}$ contains $\tau$ as a subsequence,
\textit{i.e.\/}, $1/k!$.  Thus $f(n,\tau)=1/k!$ for every
$\tau\in\S{k}$ and even $n\geq 2k$, so the limit in
Equation~\ref{eq:pqrand} equals $0$ and $\{\mathcal{P}_{n}\}$ is a
$p$-quasirandom family.
\end{ex}

As illustrated by this example, $p$-quasirandomness is weaker than
quasirandomness because the contained permutation $\tau$ is fixed
before taking a limit. The following proposition shows that the
fixed word $\tau$ need not be a permutation.
\begin{prop}
If $\{\mathcal{P}_n\}$ is a $p$-quasirandom family of
permutations, then for any fixed word
$\sigma=\sigma_1\ldots\sigma_k$ of $k$ distinct letters the
limiting probability that $\pi\in\mathcal{P}_n$ contains $\sigma$
is $1/k!$ as $n\rightarrow\infty$.
\end{prop}
\begin{proof}
Let $m$ be the largest letter which appears in the word $\sigma$.
The probability that $\pi\in\mathcal{P}_n$ contains $\pi$ is
$1/m!+o(1)$ as $n\rightarrow\infty$ for each of the
$\binom{m}{k}(m-k)! = m!/k!$ permutations which contain $\sigma$
as a subsequence.
\end{proof}

\subsection{The $p$-quasirandomness of involutions}
\label{subsec:finv}

Because the set of involutions is $p$-quasirandom, the probability
that an $n$-involution contains a given $k$-permutation as a
subsequence is $1/k!+o(1)$ as $n\rightarrow\infty$.  We now
sharpen this result by obtaining an exact count of the
$n$-involutions which contain a given $k$-permutation as a
subsequence, starting with the following lemma.
\begin{lemma}
For $\tau=\tau_{1}\ldots\tau_{k}\in\S{k}$ and $0\leq j\leq k\leq
n$, the number of $n$-involutions $\pi$ which contain $\tau$ as a
subsequence and which map exactly $j$ elements of $[k]$ into $[k]$
is $\binom{n-k}{k-j}t_{n-2k+j}$, if either $j=0$ or the pattern of
$\tau_{1}\ldots\tau_{j}$ is a $j$-involution, and $0$ otherwise.
\label{lem:invj}
\end{lemma}
\begin{proof}
Fix $\tau\in\S{k}$ and an $n$-involution $\pi$ containing $\tau$
as a subsequence.  Let $A=\{a_{1},\ldots,a_{j}\}$ be the elements
of $[k]$ which are mapped by $\pi$ into $[k]$, with
$a_{1}<a_{2}<\cdots<a_{j}$ if $A\neq\emptyset$.

Now assume that $j\neq 0$.  Since $\pi(a_{i})\in[k]$ by definition
and $\pi(\pi(a_{i})) = a_{i}\in[k]$, $\pi(a_{i})\in A$ and the
restriction of $\pi$ to $A$ is an involution in the group of
permutations of $A$. Because $\pi$ contains $\tau$ as a
subsequence and $a_{1}$ is the smallest element of $[k]$ which is
mapped by $\pi$ into $[k] = \{\tau_{1}, \ldots,\tau_{k}\}$, we
have $\pi(a_{1}) = \tau_{1}$; in general, $\pi(a_{i})=\tau_{i}\in
A$ for each $i\in[j]$.  Combining this with the fact that
$a_{1}<a_{2}<\cdots<a_{j}$ is the ordering of
$A=\{\tau_{1},\ldots,\tau_{j}\}$ in increasing order and the fact
that the restriction of $\pi$ to $A$ is an involution shows that
the pattern $x_{1}{\ldots}x_{j}$ of $\tau_{1}\ldots\tau_{j}$ is an
involution in $\S{j}$.

If $j=0$, or if $1\leq j\leq k$ and the length $j$ initial pattern
of $\tau$ is a $j$-involution, we may construct the permutations
$\pi$ which contain $\tau$ as a subsequence and which map $j$
elements of $[k]$ into $[k]$.  The requirement
$\pi(a_{i})=\tau_{i}$ noted above defines $\pi$ on $A$.  For
$\tau_{i}\notin A$, $\pi^{-1}(\tau_{i}) = \pi(\tau_{i})\notin
[k]$, so for the preimages under $\pi$ of the elements of
$[k]\setminus A$ we choose $k-j$ elements $b_{1}<\cdots<b_{k-j}$
from $[n]\setminus[k]$.  Because $\pi$ must contain the $\tau_{i}$
in the order $\tau_{1},\ldots,\tau_{k}$, we have
$\pi(b_{i})=\tau_{i+j}$, with the involuting nature of $\pi$
forcing $\pi(\tau_{i+j})=b_{i}$.  This completes the definition of
$\pi$ on the $2k-j$ elements of $[k]\cup\{b_i\}_{i\in[k-j]}$.  We
may define $\pi$ on the remaining elements of $[n]$ by choosing
one of the $t_{n-2k+j}$ involutions in the group of permutations
of $[n]\setminus\{1,\ldots,k,b_1,\ldots,b_{k-j}\}$.
\end{proof}

{From} this lemma, we may immediately count the $n$-involutions
which contain a fixed $k$-permutation $\tau$.
\begin{theorem}
For a fixed permutation $\tau=\tau_{1}\tau_{2}\ldots\tau_{k}
\in\S{k}$ and $n\geq k$, the number $\siav{n}{\tau}$ of
involutions in $\S{n}$ which contain $\tau$ as a subsequence
equals
\begin{equation}
\sideset{}{'}\sum_{j}\binom{n-k}{k-j}t_{n-2k+j} \label{eq:numinv}
\end{equation}
where the sum is taken over $j=0$ and those $j\in[k]$ such that
the pattern of $\tau_1\ldots\tau_j$ is an involution in
$\mathcal{S}_j$.\qed \label{thm:inv}
\end{theorem}
\noindent  We note that Equation~\ref{eq:numinv} is an asymptotic
series.  Using the appropriate asymptotic expansions, we may use
this theorem to give another proof of the $p$-quasirandomness of
the set of involutions and to sharpen the associated limiting
probability.
\begin{cor}
The probability that an $n$-involution contains $\tau\in\S{k}$ as
a subsequence equals
\begin{equation}
\sideset{}{'}\sum_{j}\binom{n-k}{k-j} \frac{t_{n-2k+j}}{t_{n}},
\label{eq:finv}
\end{equation}
where the sum is taken over $j=0$ and those $j\in[k]$ such that
the pattern of $\tau_1\ldots\tau_j$ is an involution in
$\mathcal{S}_j$.\qed
\end{cor}

\begin{remark}
For every $\tau$ in $\S{k}$, $k\geq 2$, the sum in
Equation~\ref{eq:finv} is taken over at least $j=0$, $1$, and
$2$.\label{rem:jsum}
\end{remark}
\begin{cor}
The probability that an $n$-involution contains the subsequence
$1$ equals $1$ for every positive value of $n$.  For $n\geq 2$ and
$\tau\in\S{2}$, the probability that an $n$-involution contains
$\tau$ as a subsequence is exactly $1/2$. \label{cor:invpqrand12}
\end{cor}

\begin{cor}
For $k>2$ and $\tau\in\S{k}$, a lower bound for the probability
that an $n$-involution contains $\tau$ as a subsequence is given
by
\begin{equation}
\binom{n-k}{k}\frac{t_{n-2k}}{t_{n}} +
\binom{n-k}{k-1}\frac{t_{n-2k+1}}{t_{n}}+
\binom{n-k}{k-2}\frac{t_{n-2k+2}}{t_{n}} \leq f_{inv}(n,\tau).
\label{eq:prob:lb}
\end{equation}
Furthermore, it is possible for equality to hold.
\end{cor}
\begin{proof}
The bound follows from Remark~\ref{rem:jsum} above.  The
permutation $k12\ldots(k-1)$ has length $j$ initial pattern
$j12\ldots(j-1)$, which is not an involution in $\S{j}$ for $j>2$;
equality holds for this permutation.
\end{proof}

\begin{cor}
For $k>2$ and $\tau\in\S{k}$, an upper bound for the probability
that an $n$-involution contains $\tau$ as a subsequence is given
by
\begin{equation}
f_{inv}(n,\tau) \leq \sum_{j=0}^k\binom{n-k}{k-j}
\frac{t_{n-2k+j}}{t_{n}}.\label{eq:prob:ub}
\end{equation}
Furthermore, it is possible for equality to hold.
\end{cor}
\begin{proof}
The bound follows immediately from Equation~\ref{eq:finv}, with
equality holding for, \textit{e.g.\/}, $\tau=12\ldots k$.
\end{proof}
\begin{remark}
Proposition~\ref{prop:c:fulljset} below shows that there are
exactly $2^{k-1}$ permutations $\tau\in\S{k}$ for which equality
holds in Equation~\ref{eq:prob:ub}.
\end{remark}

Using the appropriate asymptotic expansions, we may obtain a
sharper value than $1/k!$ for the limiting probability that an
$n$-involution contains $\tau\in\S{k}$.

\begin{prop}
For $k>2$, $\tau\in\S{k}$, the probability as $n\rightarrow\infty$
that an $n$-involution $\pi$ contains $\tau$ as a subsequence is
\begin{equation}
\frac{1}{k!}-\frac{2}{3(k-3)!}n^{-3/2}+O(n^{-2})
\end{equation}
if the pattern of $\tau_1\tau_2\tau_3$ is not an involution in
$\S{3}$ and
\begin{equation}
\frac{1}{k!}+\frac{1}{3(k-3)!}n^{-3/2}+O(n^{-2})
\end{equation}
if the pattern of $\tau_1\tau_2\tau_3$ is an involution in
$\S{3}$. \label{prop:finv}
\end{prop}

\begin{remark}
Theorem~\ref{thm:inv} and its corollaries have natural analogues
for fixed-point-free involutions~\cite{adjthesis}.
\end{remark}

\section{Applications}
\label{sec:app}

\subsection{Classifying permutations}

We now consider when permutations are equally restrictive with
respect to subsequence containment by involutions.  This has
strong parallels to the notion of \textit{Wilf-equivalence} from
the study of pattern-avoiding permutations.
\begin{defn}
Say that two permutations $\sigma_1$ and $\sigma_2$ are
\emph{equivalent with respect to subsequence containment by
involutions} iff for every $n$, the number of $n$-involutions
which contain $\sigma_1$ as a subsequence equals the number which
contain $\sigma_2$ as a subsequence. \label{def:subseqequiv}
\end{defn}
\noindent Note that the analogous definition for containment by
permutations leads to a trivial equivalence.

Our classification of permutations using this equivalence will
make use of Proposition~\ref{prop:equivsubseq} below, for which we
need the following definition.
\begin{defn}
For $\tau\in\S{k}$, the \textit{$j$-set of $\tau$}, denoted
$\J(\tau)$, is the set containing $0$ and exactly those $i\in[k]$
such that the length $i$ initial pattern of $\sigma$ is an
involution.  This is the set of indices $j$ over which the sum in
Equation~\ref{eq:numinv} is taken when counting the
$n$-involutions which contain $\tau$ as a subsequence.
\end{defn}

\begin{prop}
Two permutations are equivalent with respect to subsequence
containment by involutions iff they are of the same length and
their $j$-sets are identical. \label{prop:equivsubseq}
\end{prop}
\begin{proof}
Assume that for two distinct permutations $\sigma\in\S{k}$ and
$\tau\in\S{k'}$ the corresponding sums in Equation~\ref{eq:numinv}
are equal for every value of $n$.  Because the limiting value
($n\rightarrow\infty$) of these sums divided by $t_n$ equals
$1/k!$ and $1/k'!$, respectively, we must have $k'=k$.  As noted
above, these sums are asymptotic series, so the values of $j$ from
$\{0,1,\ldots,k\}$ used in each sum must be the same.
\end{proof}

Proposition~\ref{prop:equivsubseq} allows us to classify
permutations based on subsequence containment, \textit{i.e.\/}, to
determine the equivalence classes with respect to
Definition~\ref{def:subseqequiv}, by simply determining which
permutations have the same $j$-sets.  Table~\ref{tab:s3invsubseq}
lists the $2$ possible $j$-sets for permutations in $\S{3}$, the
number of permutations which have each of those $j$-sets, and the
number of $n$-involutions which contain the permutations from
these classes as subsequences for $3\leq n\leq 10$.\ \ (Recall
that $\siav{n}{\tau}$ denotes the number of $n$-involutions that
contain $\tau$ as a subsequence.)\ \ Table~\ref{tab:s4invsubseq}
does the same for permutations in $\S{4}$; in this case, there are
$4$ possible $j$-sets and the number of containing $n$-involutions
is given for $4\leq n\leq 10$.

\begin{table}
    \small
  \centering
  \begin{tabular}{|c|c|c|c|c|c|c|c|c|c|}
    \hline
    $\J(\tau)$ & $|\{\tau\}|$ & $\siav{3}{\tau}$ &
    $\siav{4}{\tau}$ & $\siav{5}{\tau}$ & $\siav{6}{\tau}$ &
    $\siav{7}{\tau}$ & $\siav{8}{\tau}$ & $\siav{9}{\tau}$ & $\siav{10}{\tau}$ \\ \hline
    $\{0,1,2\}$  & 2 & 0 & 1 & 3 & 10 & 32 & 110 & 386 & 1428 \\ \hline
    $\{0,1,2,3\}$ & 4 & 1 & 2 & 5 & 14 & 42 & 136 & 462 & 1660 \\ \hline
  \end{tabular}
  \caption{Classifying $\S{3}$ by subsequence containment by involutions.}\label{tab:s3invsubseq}
\end{table}

\begin{table}
    \small
  \centering
  \begin{tabular}{|c|c|c|c|c|c|c|c|c|}
    \hline
    $\J(\tau)$ & $|\{\tau\}|$ &
    $\siav{4}{\tau}$ & $\siav{5}{\tau}$ & $\siav{6}{\tau}$ &
    $\siav{7}{\tau}$ & $\siav{8}{\tau}$ & $\siav{9}{\tau}$ & $\siav{10}{\tau}$ \\ \hline
    $\{0,1,2\}$  & 6 & 0 & 0 & 1 & 4 & 17 & 65 & 260 \\ \hline
    $\{0,1,2,3\}$  & 8 & 0 & 1 & 3 & 10 & 33 & 115 & 416 \\ \hline
    $\{0,1,2,4\}$ & 2 & 1 & 1 & 3 & 8 & 27 & 91 & 336 \\ \hline
    $\{0,1,2,3,4\}$  & 8 & 1 & 2 & 5 & 14 & 43 & 141 & 492 \\ \hline
  \end{tabular}
  \caption{Classifying $\S{4}$ by subsequence containment by involutions.}\label{tab:s4invsubseq}
\end{table}

Tables~\ref{tab:s5invsubseq}--\ref{tab:s7invsubseq} present
similar data for permutations in $\S{5}$, $\S{6}$, and $\S{7}$.
Note that for $\S{7}$, not every possible $j$-set is realized as
$\J(\tau)$ for some $\tau$; neither $\{0,1,2,5,7\}$ nor
$\{0,1,2,4,5,7\}$ equals $\J(\tau)$ for some $\tau\in\S{7}$.
Additionally, no permutation $\tau\in\S{8}$ has
\begin{multline*}
\J(\tau)\in \{\{0,1,2,5,7\}, \{0,1,2,6,8\}, \{0,1,2,3,6,8\},
\{0,1,2,4,5,7\},\\ \{0,1,2,5,6,8\}, \{0,1,2,5,7,8\},
\{0,1,2,3,5,6,8\}, \{0,1,2,4,5,7,8\} \}.
\end{multline*}
This suggests that it may be interesting to determine how many
$j$-sets actually occur.
\begin{q}
What is the sequence
\begin{equation*}
\{|\J(\S{k})|\}_{k\geq 3} = 2,4,8,16,30,56,102,\ldots?
\end{equation*}
\textit{I.e.\/}, for $k\geq 3$, how may of the $2^{k-2}$ possible
$j$-sets are actually realized by some permutation in $\S{k}$?
\end{q}
A more general question is the enumeration of the $k$-permutations
which have a particular $j$-set (as opposed to simply determining
when this count is nonzero).
\begin{q}
Given a set $E$, $\{0,1,2\}\subseteq E\subseteq \{0,1,\ldots,k\}$,
how many $k$-permutations $\tau$ have $\J(\tau)=E$?
\end{q}

\begin{table}[htb]
  \centering
  \begin{tabular}{|c|c|c|c|c|c|c|c|c|c|}
    \hline
    $\J(\tau)$ & $|\{\tau\}|$ & $\siav{5}{\tau}$ & $\siav{6}{\tau}$ &
    $\siav{7}{\tau}$ & $\siav{8}{\tau}$ & $\siav{9}{\tau}$ & $\siav{10}{\tau}$ \\ \hline
    $\{0,1,2\}$  & 26 & 0 & 0 & 0 & 1 & 5 & 26 \\ \hline
    $\{0,1,2,3\}$ & 36 & 0 & 0 & 1 & 4 & 17 & 66 \\ \hline
    $\{0,1,2,4\}$ & 8 & 0 & 1 & 2 & 7 & 21 & 76 \\
    \hline
    $\{0,1,2,5\}$ & 4 & 1 & 1 & 2 & 5 & 15 & 52 \\
    \hline
    $\{0,1,2,3,4\}$ & 24 & 0 & 1 & 3 & 10 & 33 & 116 \\
    \hline
    $\{0,1,2,3,5\}$ & 4 & 1 & 1 & 3 & 8 & 27 & 92 \\
    \hline
    $\{0,1,2,4,5\}$ & 2 & 1 & 2 & 4 & 11 & 31 & 102 \\
    \hline
    $\{0,1,2,3,4,5\}$ & 16 & 1 & 2 & 5 & 14 & 43 & 142 \\
    \hline
  \end{tabular}
  \caption{Classifying $\S{5}$ by subsequence containment by involutions.}\label{tab:s5invsubseq}
\end{table}

\begin{table}[htb]
  \centering
  \begin{tabular}{|c|c|c|c|c|c|c|c|c|c|}
    \hline
    $\J(\tau)$ & $|\{\tau\}|$ & $\siav{6}{\tau}$ &
    $\siav{7}{\tau}$ & $\siav{8}{\tau}$ & $\siav{9}{\tau}$ & $\siav{10}{\tau}$ \\ \hline
    $\{0,1,2\}$  & 146 & 0 & 0 & 0 & 0 & 1 \\ \hline
    $\{0,1,2,3\}$ & 204 & 0 & 0 & 0 & 1 & 5
    \\ \hline
    $\{0,1,2,4\}$ & 46 & 0 & 0 & 1 & 3 & 13 \\
    \hline
    $\{0,1,2,5\}$ & 20 & 0 & 1 & 2 & 6 & 17 \\
    \hline
    $\{0,1,2,6\}$ & 10 & 1 & 1 & 2 & 4 & 11 \\
    \hline
    $\{0,1,2,3,4\}$ & 136 & 0 & 0 & 1 & 4 & 17 \\
    \hline
    $\{0,1,2,3,5\}$ & 20 & 0 & 1 & 2 & 7 & 21 \\
    \hline
    $\{0,1,2,3,6\}$ & 12 & 1 & 1 & 2 & 5 & 15 \\
    \hline
    $\{0,1,2,4,5\}$ & 8 & 0 & 1 & 3 & 9 & 29 \\
    \hline
    $\{0,1,2,4,6\}$ & 2 & 1 & 1 & 3 & 7 & 23 \\
    \hline
    $\{0,1,2,5,6\}$ & 4 & 1 & 2 & 4 & 10 & 27 \\
    \hline
    $\{0,1,2,3,4,5\}$ & 64 & 0 & 1 & 3 & 10 & 33 \\
    \hline
    $\{0,1,2,3,4,6\}$ & 8 & 1 & 1 & 3 & 8 & 27 \\
    \hline
    $\{0,1,2,3,5,6\}$ & 4 & 1 & 2 & 4 & 11 & 31 \\
    \hline
    $\{0,1,2,4,5,6\}$ & 4 & 1 & 2 & 5 & 13 & 39 \\
    \hline
    $\{0,1,2,3,4,5,6\}$ & 32 & 1 & 2 & 5 & 14 & 43 \\
    \hline
  \end{tabular}
  \caption{Classifying $\S{6}$ by subsequence containment by involutions.}\label{tab:s6invsubseq}
\end{table}

\begin{table}
  \centering
  \begin{tabular}{|c|c|c|c|c|c|c|c|c|c|}
    \hline
    $\J(\tau)$ & $|\{\tau\}|$ &
    $\siav{7}{\tau}$ & $\siav{8}{\tau}$ & $\siav{9}{\tau}$ & $\siav{10}{\tau}$ \\ \hline
    $\{0,1,2\}$  & 992 & 0 & 0 & 0 & 0 \\ \hline
    $\{0,1,2,3\}$ & 1396 & 0 & 0 & 0 & 0
    \\ \hline
    $\{0,1,2,4\}$ & 316 & 0 & 0 & 0 & 1 \\
    \hline
    $\{0,1,2,5\}$ & 140 & 0 & 0 & 1 & 3 \\
    \hline
    $\{0,1,2,6\}$ & 60 & 0 & 1 & 2 & 6 \\
    \hline
    $\{0,1,2,7\}$ & 30 & 1 & 1 & 2 & 4 \\
    \hline
    $\{0,1,2,3,4\}$ & 928 & 0 & 0 & 0 & 1 \\
    \hline
    $\{0,1,2,3,5\}$ & 136 & 0 & 0 & 1 & 3 \\
    \hline
    $\{0,1,2,3,6\}$ & 72 & 0 & 1 & 2 & 6 \\
    \hline
    $\{0,1,2,3,7\}$ & 32 & 1 & 1 & 2 & 4 \\
    \hline
    $\{0,1,2,4,5\}$ & 56 & 0 & 0 & 1 & 4 \\
    \hline
    $\{0,1,2,4,6\}$ & 12 & 0 & 1 & 2 & 7 \\
    \hline
    $\{0,1,2,4,7\}$ & 6 & 1 & 1 & 2 & 5 \\
    \hline
    $\{0,1,2,5,6\}$ & 20 & 0 & 1 & 3 & 9 \\
    \hline
    $\{0,1,2,5,7\}$ & 0 & 1 & 1 & 3 & 7 \\
    \hline
    $\{0,1,2,6,7\}$ & 10 & 1 & 2 & 4 & 10 \\
    \hline
    $\{0,1,2,3,4,5\}$ & 432 & 0 & 0 & 1 & 4 \\
    \hline
    $\{0,1,2,3,4,6\}$ & 48 & 0 & 1 & 2 & 7 \\
    \hline
    $\{0,1,2,3,4,7\}$ & 24 & 1 & 1 & 2 & 5 \\
    \hline
    $\{0,1,2,3,5,6\}$ & 20 & 0 & 1 & 3 & 9 \\
    \hline
    $\{0,1,2,3,5,7\}$ & 4 & 1 & 1 & 3 & 7 \\
    \hline
    $\{0,1,2,3,6,7\}$ & 12 & 1 & 2 & 4 & 10 \\
    \hline
    $\{0,1,2,4,5,6\}$ & 20 & 0 & 1 & 3 & 10 \\
    \hline
    $\{0,1,2,4,5,7\}$ & 0 & 1 & 1 & 3 & 8 \\
    \hline
    $\{0,1,2,4,6,7\}$ & 2 & 1 & 2 & 4 & 11 \\
    \hline
    $\{0,1,2,5,6,7\}$ & 8 & 1 & 2 & 5 & 13 \\
    \hline
    $\{0,1,2,3,4,5,6\}$ & 160 & 0 & 1 & 3 & 10 \\
    \hline
    $\{0,1,2,3,4,5,7\}$ & 16 & 1 & 1 & 3 & 8 \\
    \hline
    $\{0,1,2,3,4,6,7\}$ & 8 & 1 & 2 & 4 & 11 \\
    \hline
    $\{0,1,2,3,5,6,7\}$ & 8 & 1 & 2 & 5 & 13 \\
    \hline
    $\{0,1,2,4,5,6,7\}$ & 8 & 1 & 2 & 5 & 14 \\
    \hline
    $\{0,1,2,3,4,5,6,7\}$ & 64 & 1 & 2 & 5 & 14 \\
    \hline
  \end{tabular}
  \caption{Classifying $\S{7}$ by subsequence containment by involutions.}\label{tab:s7invsubseq}
\end{table}

As special cases of this question, we have the following question
and proposition.  These are of particular interest because they
give the number of $k$-permutations $\tau$ which achieve (for
large enough $n$) the smallest and largest values of
$\siav{n}{\tau}$.

\begin{q}
How may $k$-permutations have $j$-set equal to $\{0,1,2\}$?
\end{q}

\begin{figure}[htb]
\begin{center}
\includegraphics[width=4in]{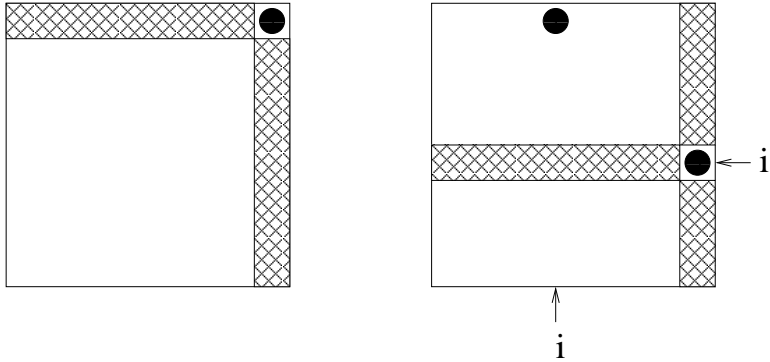}
\caption{Constructing the graphs of $\opn{\tau}$, left, and
$\opi{\tau}$, right, from the graph of $\tau$.}\label{fig:optau}
\end{center}
\end{figure}

\begin{prop}
The number of $\tau\in\S{k}$ for which $\J(\tau)=\{0,1,\ldots,k\}$
equals $2^{k-1}$. \label{prop:c:fulljset}
\end{prop}
Before proving this proposition, we define two operations for
extending a permutation $\tau\in\S{k}$ to a permutation in
$\S{k+1}$ whose length $k$ initial pattern equals $\tau$.  It will
be helpful to use the \emph{graph of a permutation} in doing so;
the graph of $\tau\in\S{k}$ is a $k\times k$ grid, which we will
co-ordinatize from the bottom left corner, with dots in exactly
the boxes $\{(i,\tau_i)\}_{i\in[k]}$.
\begin{defn}
Given $\tau\in\S{k}$, $\opn{\tau}$ is the permutation in $\S{k+1}$
that fixes $k+1$ and permutes $[k]$ as $\tau$ does.  The graph of
$\opn{\tau}$ is obtained from the graph of $\tau$ by adding a dot
in the box $(k+1,k+1)$ (as well as the appropriate additional
empty boxes).  The left part of Figure~\ref{fig:optau} illustrates
the construction of the graph of $\opn{\tau}$ from the graph of
$\tau$; the white area indicates the boxes forming the graph of
$\tau$ while the shaded areas and dot are added to obtain the
graph of $\opn{\tau}$.\label{def:opn}
\end{defn}
\begin{remark}
By the construction of $\opn{\tau}$, we see that for
$\tau\in\S{k}$:
\begin{enumerate}
    \item  the length $k$ initial pattern of $\opn{\tau}$ equals
    $\tau$ and
    \item  $\opn{\tau}$ is an involution iff $\tau$ is an involution.
\end{enumerate}
\end{remark}
\begin{defn}
Given $\tau\in\S{k}$ with $\tau_i=k$, $\opi{\tau}$ equals the
permutation obtained by adding $1$ to every value in $\tau$ that
is at least $i$ and then appending the value $i$.  The graph of
$\opi{\tau}$ is obtained from the graph of $\tau$ by inserting an
empty row at height $i$ (moving the rows originally at heights
$i,\ldots,k$ to be at heights $i+1,\ldots,k+1$) and then adding a
dot in the box $(k+1,i)$ (as well as the appropriate additional
empty boxes).  The right part of Figure~\ref{fig:optau}
illustrates the construction of the graph of $\opi{\tau}$ from the
graph of $\tau$ when $\tau_i=k$ is the maximum value of
$\tau\in\S{k}$. The white area again indicates the boxes forming
the graph of $\tau$, now split into the bottom $i-1$ rows and the
remaining $k-i+1$ rows, while the shaded boxes and dot are added
to obtain the graph of $\opi{\tau}$. \label{def:opi}
\end{defn}
\begin{remark}
By the construction of $\opi{\tau}$, we see that for
$\tau\in\S{k}$:
\begin{enumerate}
    \item   the length $k$ initial pattern of $\opi{\tau}$ equals
    $\tau$ and
    \item   $\opi{\tau}$ is an involution iff the length $k-1$
    pattern of $\tau^{-1}$ is an involution (equivalently, iff the
    permutation obtained by deleting the largest value $k$ from
    $\tau$ is an involution).
\end{enumerate}
\end{remark}
Note that if an involution in $\S{k+1}$ has length $k$ initial
pattern equal to $\tau\in\S{k}$, then that involution must equal
$\opn{\tau}$ (if it fixes $k+1$) or $\opi{\tau}$ (if it does not).
\begin{proof}
\emph{(Of Proposition~\ref{prop:c:fulljset})}\ \ This is true for
$k=1,2$; assume that the proposition holds for $k$ and pick
$\tau\in\S{k}$ such that $\J(\tau)=\{0,\ldots,k\}$. As $\tau$ is
an involution, $\opn{\tau}$ is as well and
$\J(\opn{\tau})=\{0,\ldots,k,k+1\}$. The length $k-1$ initial
pattern of $\tau^{-1}=\tau$ is also an involution by assumption,
so $\opi{\tau}$ is an involution and
$\J(\opi{\tau})=\{0,\ldots,k,k+1\}$.  Because
$\opn{\tau}\neq\opi{\tau}$ and any involution in $\S{k+1}$ whose
length $k$ initial pattern equals $\tau$ must be either
$\opn{\tau}$ or $\opi{\tau}$, the proposition holds for $k+1$.
\end{proof}

\begin{defn}
A permutation $\tau\in\S{k}$ is said to be \emph{layered} if, for
some composition $(a_1,a_2,\ldots,a_j)$ of $k$ (with each $a_i\geq
1$), $\tau$ consists of the first $a_1$ positive integers arranged
in decreasing order, then the next $a_2$ positive integers in
decreasing order, \emph{etc.}  For example, $32146578\in\S{8}$ is
the layered permutation corresponding to the composition
$(3,1,2,1,1)$ of $8$.
\end{defn}

\begin{prop}
A permutation $\tau\in\S{k}$ is a layered permutation if and only
if $\J(\tau)=\{0,1,\ldots,k\}$, \emph{i.e.}, if and only if the
pattern of $\tau_1\ldots\tau_j$ is an involution for every $j\in
[k]$.\label{prop:layered}
\end{prop}
\begin{proof}
If $\sigma\in\S{k'}$ is layered, let the last layer of $\sigma$ be
$\sigma_i=k,\ldots,\sigma_k=i$. $\opi{\sigma}$ agrees with
$\sigma$ in its first $i-1$ values; these are followed in
$\opi{\sigma}$ by $(k+1)\ldots(i+1)i$ and so $\opi{\sigma}$ is
also layered.  It is clear that if $\sigma\in\S{k'}$ is layered
then $\opn{\sigma}$ is also layered.  The proof of
Proposition~\ref{prop:c:fulljset} shows that every $\tau\in\S{k}$
with $\J(\tau)=\{0,\ldots,k\}$ may be constructed by starting with
a permutation in $\S{2}$, both of which are layered, and applying
any sequence of operations $\opn{\cdot}$ and $\opi{\cdot}$, so all
$2^{k-1}$ of these $\tau\in\S{k}$ are layered. There are $2^{k-1}$
layered permutations in $\S{k}$, a well-known result that follows
from the natural bijection between layered permutations and
compositions of $k$.
\end{proof}

\begin{lemma}
Given a layered permutation $\tau\in\S{k}$, $k\geq 2$, there is a
unique involution $\itau\in\S{k+2}$ whose length $k$ initial
pattern equals $\tau$ and whose length $k+1$ initial pattern is
not an involution.\label{lem:Ekplus2}
\end{lemma}
\begin{proof}
If such an involution $\itau$ exists, then it must equal
$\opi{\sigma}$, where $\sigma$ is the length $k+1$ initial pattern
of $\itau$; $\itau$ cannot equal $\opn{\sigma}$ as $\opn{\cdot}$
preserves the (non-)involutive nature of permutations.
$\opi{\sigma}$ is an involution iff the permutation obtained from
$\sigma$ by deleting its largest value is an involution. If
$\alpha\in\S{k}$ is obtained from $\sigma$ by deleting its largest
value, then the length $k-1$ initial pattern of $\alpha$ must be
layered.  Thus if $\alpha$ is an involution, it must be obtained
by applying $\opn{\cdot}$ or $\opi{\cdot}$ to a layered
permutation.

Note that the largest value of $\sigma$ may not be $\sigma_{k+1}$,
otherwise $\sigma$ would be layered.  If the largest value of
$\sigma$ is $\sigma_{k}$ (\emph{i.e.}, if the last layer of $\tau$
has size $1$), then $\alpha$ could not be the result of applying
$\opn{\cdot}$ to a layered permutation as $\sigma$ would then be a
layered permutation.  However, in this case we may take
$\alpha=\opi{\tau_1\ldots\tau_{k-1}}$ and obtain $\sigma$ by
inserting $k+1$ just before the last value of $\alpha$;
$\itau=\opi{\sigma}$ is then as required.

If the largest value of $\sigma$ is not $\sigma_k$ (\emph{i.e.},
if the last layer of $\tau$ has size greater than $1$), then
$\alpha$ could not be the result of applying $\opi{\cdot}$ to a
layered permutation as $\sigma$ would then be a layered
permutation.  However, if $j$ is then index of the largest value
of $\tau$, we may take
$\alpha=\opn{\tau_1\ldots\tau_{j-1}\tau_{j+1}\ldots\tau_k}$ and
obtain $\sigma$ by inserting $k+1$ into $\alpha$ just before
$\alpha_j$; $\itau=\opi{\sigma}$ is then as required.

As $\itau$ must equal $\opi{\sigma}$ and there is a unique way to
extend any layered $\tau\in\S{k}$ to $\sigma$ such that $\itau$
has the stated properties (the manner of doing so depending upon
whether the last layer of $\tau$ has size $1$), the lemma is
proved.
\end{proof}

\begin{lemma}
Given $\tau\in\S{k}$ with $\J(\tau)=\{0,1,\ldots,k\}$, there is a
unique $\jtau\in\S{k+3}$ with $\J(\jtau)=\{0,1,\ldots,k,k+2,k+3\}$
and length $k$ initial pattern equal to $\tau$.\label{lem:Ekplus3}
\end{lemma}
\begin{proof}
Construct the $\itau\in\S{k+2}$ as in Lemma~\ref{lem:Ekplus2};
this is the only possible length $k+2$ initial pattern of $\jtau$.
As $\itau$ is an involution, $\opn{\itau}$ is as well. By
construction, the length $k+1$ initial pattern of
$\itau^{-1}=\itau$ is not an involution, so $\opi{\itau}$ is not
an involution.  $\opi{\itau}$ is the only other way that $\itau$
could have been extended to an involution in $\S{k+3}$, so
$\jtau=\opn{\itau}$ is unique as claimed.
\end{proof}

Tables~\ref{tab:s3invsubseq}--\ref{tab:s7invsubseq} suggest the
following proposition, which we may prove using these lemmas.
\begin{prop}
Assume $\{0,1,2,k\}\subseteq E\subset \{0,1,\ldots,k\}$ and
$|E|=k\geq 5$.  Then
\begin{equation*}
|j^{-1}(E)\bigcap\S{k}|=
\begin{cases}
2^{k-3} & k-1\notin E\\
2^{k-4} & k-1\in E
\end{cases}.
\end{equation*}
\end{prop}
\begin{proof}
Let $\{i\}=\{0,1,\ldots,k\}\setminus E$.  The case $i=k-1$ is
covered by Proposition~\ref{prop:c:fulljset} and
Lemma~\ref{lem:Ekplus2}, while the case $i=k-2$ is covered by
these and Lemma~\ref{lem:Ekplus3}.

For $i\leq k-3$, construct the unique $\jtau\in\S{i+2}$ such that
$\J(\jtau)=\{0,\ldots,i-1,i+1,i+2\}$.  Because $\jtau$ is an
involution, the length $i+1$ initial pattern of $\jtau^{-1}$
equals the length $i+1$ initial pattern of $\jtau$ itself; by
assumption, this is an involution.  Thus both $\opn{\jtau}$ and
$\opi{\jtau}$ are involutions in $\S{i+3}$ whose $j$-sets equal
$\{0,\ldots,i-1,i+1,i+2,i+3\}$, and the proposition holds for
$k=i+3$.  The same argument shows that these permutations may be
extended by any sequence of $\opn{\cdot}$ and $\opi{\cdot}$
operations until length $k$ involutions are produced with the
desired properties, and the theorem is proved.
\end{proof}

Finally, review of
Tables~\ref{tab:s3invsubseq}--\ref{tab:s7invsubseq} also suggests
the following question.
\begin{q}
For $k>4$ and any set $E\neq\{0,1,2,3\}$, is the number of
$\sigma\in\S{k}$ for which $\J(\sigma)=E$ always strictly less
than the number of $\sigma\in\S{k}$ for which
$\J(\sigma)=\{0,1,2,3\}$?  Whether or not this holds, are there
any other statements that can be made about the frequency of
particular sets appearing as $\J(\sigma)$?
\end{q}

\subsection{Tableau containment}

We might ask similar questions for certain classes of permutations
instead of for all permutations in $\S{n}$.  As an example, fix a
(standard) Young tableau $T$ of size $k$ and consider the set $\T$
of $k$-permutations which correspond under the
Robinson--Schensted-Knuth (RSK) algorithm to $(T,Q)$ for some
tableau $Q$ of the same shape as $T$. Lemma~\ref{lem:aj} counts
the permutations $\sigma\in\T$ such that $\J(\sigma)$ contains a
specified value; we first briefly recall the RSK algorithm and
some of its properties.

The RSK algorithm bijectively associates to every
$n$-per\-mu\-tation $\pi$ a pair $(P,Q)$ of standard Young
tableaux, where $P$ and $Q$ have as their common shape some
partition of $n$.  Following~\cite{EC2}, we refer to $P$ as the
\textit{insertion  tableau} and to $Q$ as the \textit{recording
tableau}.  In the case that  $\pi$ is a general word of $n$
distinct letters, the RSK algorithm  produces a pair $(P,Q)$,
where $P$ is a tableau whose entries are the  letters appearing in
$\pi$ and $Q$ is a standard Young tableau of the same shape as
$P$.  Applying the RSK algorithm to the pattern of $\pi$ yields
$(P',Q')$, where $P'$ is the order preserving relabeling of $P$
with the elements of $[n]$ and $Q'=Q$. Finally, for a permutation
$\pi\in\S{n}$ with corresponding pair of tableaux $(P,Q)$, $P=Q$
iff $\pi$ is an involution.

\begin{lemma}
Fix $k$ and $\lambda\vdash k$, and let $T$ be a standard Young
tableau of shape $\lambda$.  Let $\T$ be the set of permutations
which correspond to $(T,Q)$ for some $Q$, let $a_{0}=f^{\lambda}$,
and for $1\leq j\leq k$ let $a_{j}$ be the number of
$\tau\in\mathcal{T}$ for which $j\in\J(\tau)$.  Then (with
$f^{\emptyset}=1$)
\begin{equation}
a_{j}=\sum_{\mu\,\vdash\, j}f^{\lambda/\mu} \quad (j=0,1,\ldots,
k),
\end{equation}
where $f^{\lambda/\mu}$ is the number of standard Young tableaux
of skew shape $\lambda/\mu$. \label{lem:aj}
\end{lemma}
\begin{proof}
This is immediate for $j=0$.  For $j\in[k]$ and a tableau $Q$ of
size $k$, let $\Qj$ denote the subtableau of $Q$ formed by the
elements of $[j]$.  Observe that for $j\in[k]$, if
$\tau=\tau_1\ldots\tau_k$ is a $k$-permutation which corresponds
under RSK to the pair $(P,Q)$ of tableaux, then the following
three conditions are equivalent:
\begin{enumerate}
\item[(i)] The length $j$ initial pattern of $\tau$ is a
$j$-involution.

\item[(ii)] The length $j$ initial pattern of $\tau$ corresponds
under RSK to $(\Qj,\Qj)$.

\item[(iii)] $\tau_{1}\ldots \tau_{j}$ corresponds to $(P',Q')$
under RSK, where $Q'=\Qj$ is the order preserving relabeling of
$P'$ with the elements of $[j]$.
\end{enumerate}

Each of the $f^{\lambda}$ permutations in $\mathcal{T}$
corresponds under RSK to $(T,Q)$ for some recording tableau $Q$ of
shape $\lambda$.  There are $\sum_{\mu\vdash  j}f^{\lambda/\mu}$
different ways that a recording tableau can contain
$\{j+1,\ldots,k\}$.  For a fixed such arrangement $\st$, a skew
tableau, let $\stc$ denote the set of (recording) tableaux which
contain $\{j+1,\ldots,k\}$ in the arrangement $\st$.  If $\mu$ is
the partition for which the shape of $\st$ is $\lambda/\mu$, then
there are $f^{\mu}$ tableaux in $\stc$. For a fixed set $\stc$ and
$Q\in\stc$, the removal of $k-j$ elements from $T$  by applying
the inverse RSK algorithm to $(T,Q)$ yields a pair $(T_{\st}',Q')$
of tableaux of shape $\mu$, where $Q'=\Qj$  and $T_{\st}'$ depends
only on $\st$.  Of the $f^{\mu}$ tableaux $Q'$ so obtained by
letting $Q$ range over $\st$, exactly one is the order preserving
relabeling of $T_{\st}'$ with the elements of $[j]$.  By the
observation above, this choice of $Q'$, together with  the
arrangement $\st$, gives a recording tableaux $Q\in\stc$ such that
$(T,Q)$ corresponds under RSK to a permutation whose length  $j$
initial pattern is a $j$-involution; this is the only choice of
$Q\in\stc$ for which this is true. Considering all possible
arrangements $\st$ of $\{j+1,\ldots,k\}$, we see that there are
$\sum_{\mu\vdash j}f^{\lambda/\mu}$ $k$-permutations in
$\mathcal{T}$  whose length $j$ initial pattern in a
$j$-involution.
\end{proof}

\begin{remark}
The values of $a_{j}$ in the lemma depend only on the shape
$\lambda$ of $T$.  However, the number of permutations whose
length $j$ and $j'\neq j$ initial patterns are both involutions
need not be the same for different tableaux of the same shape.
\end{remark}

As a corollary, this allows us to count the standard Young
tableaux which contain a given tableau.  Theorem~\ref{thm:tab} was
first proved by Sagan and Stanley, following from Corollary~3.5
in~\cite{Sagan1990}. The position of certain elements in a tableau
was the motivating question in the work of McKay, Morse, and Wilf
on quasirandomness~\cite{McKay2000}.  More recently, Stanley
revisited the subtableau question and developed related asymptotic
work~\cite{Stanley2001}.  Finally, Grabiner has used a random-walk
approach to study the asymptotics for the problem and a
generalization to up-down tableaux~\cite{Grabiner2004}.

\begin{theorem}
Let $T$ be a standard Young tableau of shape $\lambda\vdash k$.
Then the number of tableaux of size $n\geq k$ which contain $T$ as
a subtableau equals
\begin{equation}
\sum_{j=0}^{k} \sum_{\mu\,\vdash\, j}f^{\lambda/\mu}
\binom{n-k}{k-j}t_{n-2k+j}.
\end{equation}
\label{thm:tab}
\end{theorem}
\begin{proof}
The number of tableaux containing $T$ equals the sum of
Equation~\ref{eq:numinv} over all permutations $\tau$ which
correspond to $(T,Q)$ for some $Q$. For $j=0,1,\ldots,k$, there
are exactly $a_{j}=\sum_{\mu\,\vdash\, j}f^{\lambda/\mu}$
permutations $\tau\in\mathcal{T}$ for which the range  of the sum
in Equation~\ref{eq:numinv} includes $j$, from which the theorem
follows.
\end{proof}

\section*{Acknowledgments}

We are grateful to Herb Wilf for suggesting this problem to us and
for helpful discussions, the referee for information about related
work and other suggestions, Alkes Price for comments about known
results on layered patterns, and Andre Scedrov for arranging
research support.

\providecommand{\bysame}{\leavevmode\hbox
to3em{\hrulefill}\thinspace}
\providecommand{\MR}{\relax\ifhmode\unskip\space\fi MR }
\providecommand{\MRhref}[2]{%
  \href{http://www.ams.org/mathscinet-getitem?mr=#1}{#2}
} \providecommand{\href}[2]{#2}


\begin{thebibliography}{MMW02}

\bibitem[Gra04]{Grabiner2004}
David Grabiner, \emph{Asymptotics for the distributions of
subtableaux in
  {Y}oung and up-down tableaux}, Preprint, 2004.

\bibitem[Jag03]{adjthesis}
Aaron~D. Jaggard, \emph{Involutions in the symmetric group:
Containment
  properties and parallels to general permutations}, Ph.D. thesis, University
  of Pennsylvania, 2003.

\bibitem[MMW02]{McKay2000}
Brendan~D. McKay, Jennifer Morse, and Herbert~S. Wilf, \emph{The
distributions
  of the entries of {Y}oung tableaux}, J. Combin. Theory Ser. A \textbf{97}
  (2002), no.~1, 117--128.

\bibitem[SS90]{Sagan1990}
Bruce~E. Sagan and Richard~P. Stanley, \emph{Robinson-{S}chensted
algorithms
  for skew tableaux}, J. Combin. Theory Ser. A \textbf{55} (1990), no.~2,
  161--193.

\bibitem[Sta99]{EC2}
Richard~P. Stanley, \emph{Enumerative combinatorics. {V}ol. 2},
Cambridge
  University Press, Cambridge, 1999.

\bibitem[Sta03]{Stanley2001}
\bysame, \emph{On the enumeration of skew {Y}oung tableaux}, {Adv.
in Appl.
  Math.} \textbf{30} (2003), 283--294.

\end{thebibliography}
\end{document}